\documentclass[twocolumn,showpacs,preprintnumbers,amsmath,amssymb]{revtex4}

\usepackage{amsmath,amssymb}
\usepackage{graphicx}
\newcommand{\be}{\begin{equation}}
\newcommand{\ee}{\end{equation}}
\newcommand{\mk}{\langle k\rangle}
\newcommand{\bk}{{\bf k}}

\begin{document}
\title[]
      {The effects of degree distributions in random networks of Type-I neurons}
\author{Carlo R. Laing}
\address{School of Natural and Computational Sciences\\
         Massey University\\
         Private Bag 102-904 NSMC, Auckland, New Zealand} 
\email{c.r.laing@massey.ac.nz}
\date{\today}
\begin{abstract}
We consider large networks of theta neurons and use the Ott/Antonsen ansatz to derive
degree-based mean field equations governing the expected dynamics of the networks. Assuming
random connectivity we investigate the effects of varying the widths of the
in- and out-degree distributions on the dynamics of excitatory or inhibitory synaptically
coupled networks, and gap junction coupled networks. For synaptically coupled networks,
the dynamics are independent of the out-degree distribution. Broadening the
in-degree distribution destroys oscillations in inhibitory networks and decreases
the range of bistability in excitatory networks. For gap junction coupled neurons, broadening
the degree distribution varies the values of parameters at which 
there is an onset of collective oscillations. Many of the results are shown to also occur
in networks of more realistic neurons.

\end{abstract}


\maketitle

\section{Introduction}

It is well-known that the structure of a network can have a significant effect on its
dynamics. One type of network of great interest is that of networks of 
neurons, and much effort has gone into investigating this 
issue~\cite{rox11,schkih15,nykfri17,marhou16,perdeg13}.
In networks of neurons, the nodes are normally thought of an individual neurons, while the
edges describe connections between neurons. These connections can be directed
(in the case of synaptic connections~\cite{ermter10}) or undirected (in the case of gap junction 
connections~\cite{Benzuk04,erm06}).
One of the important properties of a node in a network is its degree: a node's in-degree is the
number of connections {\it to} it, and a node's out-degree is the number of connections {\it from}
it. In the case of undirected connections a node simply has a degree, as there is no
distinction between incoming or outgoing connections. Previous work on understanding
the dynamics of networks of
neurons has considered
the effects of correlations between the in- and out-degrees of individual
neurons~\cite{laibla20,lamsmi10,marhou16,vashou12,vegrox19,nykfri17} 
and of degree assortativity, in which the probability that
two neurons are connected is influenced by the degrees of the two 
neurons~\cite{blasche2020,kahsok17,schkih15,defranciscis2011}.

In this paper we consider the effects of varying the widths of
the distributions of in-degrees and out-degrees
on large, randomly connected networks of Type-I neurons, i.e.,~neurons whose onset of firing
is through a saddle-node-on-invariant-circle (SNIC) bifurcation. The analysis is done of
networks of theta neurons, as these are the normal form of the SNIC bifurcation. They also
have the property of being amenable to the use of the Ott/Antonsen ansatz~\cite{ottant08,ottant09},
a now well-used method for deriving equations governing the evolution of order parameter-like
quantities, valid for large networks with particular forms of 
heterogeneity~\cite{lai14A,lai09A,lukbar13,resott14,chahat17}.

A common assumption in creating random networks of neurons is that there is a fixed probability
of connecting any two neurons~\cite{tatolm12,golde10}. 
Such networks are often referred to as Erd{\"os}-R{\'e}nyi
and have a binomial distribution of degrees, with the ratio of mean degree
to standard deviation of degrees going to zero as network size goes to infinity.
However, real networks of neurons are observed to have properties incompatible
with this assumption~\cite{sonsjo05}.

Previous work on investigating the effects of degree distributions includes that of
Roxin~\cite{rox11}. He found in an inhibitory network of leaky integrate-and-fire
neurons that
broadening the distribution of in-degrees suppressed macroscopic oscillations. 
This result was reproduced in a simplified rate model which included the heterogeneity
in neuronal input due to the in-degree of cells. He also considered a network with both
excitatory and inhibitory neurons and investigated the effects of varying the
in- and out-degree distributions of the recurrent excitatory connections.
In other work,~\cite{qinche18} studied the effects of degree distribution in feedforward networks.

We consider synaptic coupling in Sec.~\ref{sec:syn} and gap junction coupling in
Sec.~\ref{sec:gap}. Both sections start with a derivation of the relevant equations using the
Ott/Antonsen ansatz and then give some numerical results. Most results use Lorentzian
distributions of heterogeneous parameters, and in Sec.~\ref{sec:gauss} we briefly
discuss results using normal and uniform distributions. We conclude in Sec.~\ref{sec:summ}.
Appendix~\ref{sec:app} contains the description of the network
of Morris-Lecar neurons used to verify some of the results derived for networks of theta
neurons.

\section{Synaptic coupling}
\label{sec:syn}
We first consider networks of theta neurons coupled by synaptic input currents with a timescale
$\tau$. Alternative formulations could model synaptic input by using dynamic 
synapses~\cite{coobyr19} or input current pulses with both a rise time and a decay 
time~\cite{ermter10}.

\subsection{Model and Theory}
We consider a network of $N$ theta neurons, each of which has dynamics described by
\be
   \frac{d\theta_i}{dt}=1-\cos{\theta_i}+(1+\cos{\theta_i})(\eta_i+I_i(t)) \label{eq:dthetadt}
\ee
for $i=1,\dots N$,
where the input current to neuron $i$ is
\be
   I_i(t)=\frac{K}{\mk}\sum_{j=1}^N A_{ij}u_j(t)
\ee
and the synaptic variables have dynamics given by
\be
   \tau\frac{du_j}{dt}=\sum_{m\in\mathbb{Z}}\delta(t-T_j^m)-u_j, \label{eq:dudt}
\ee
where $T_j^m$ is the $m$th firing time of neuron $j$, defined to happen every time
$\theta_j$ increases through $\pi$, and $\delta(\cdot)$ is the Dirac delta. 
Thus every time $\theta_j$ increases through $\pi$, $u_j$ is instantaneously 
incremented by an amount $1/\tau$, and between firing times it decays as $\sim e^{-t/\tau}$.
$K$ is the strength of connections between neurons (which may be positive or negative),
$\mk$ is the average in-degree in the network, and $A$ describes the connectivity of the
network, i.e.,~$A_{ij}=1$ if there is a connection from neuron $j$ to neuron
$i$ and $A_{ij}=0$ otherwise. We have $\mk=\sum_{i,j}A_{i,j}/N$.
 The $\eta_i$ are randomly chosen from a Lorentzian
\be
   g(\eta)=\frac{\Delta/\pi}{(\eta-\eta_0)^2+\Delta^2}
\ee
with centre $\eta_0$ and half-width-at-half-maximum (HWHM) $\Delta$, which introduces heterogeneity
to the network. The Lorentzian is chosen so that analytical progress can be made. We
discuss results for other distributions below.

If the input current $I_i$ is constant, the theta neuron shows one of two types of behaviour.
For $\eta_i+I_i<0$,~\eqref{eq:dthetadt} has two fixed points, one stable and one unstable.
For $\eta_i+I_i>0$,~\eqref{eq:dthetadt} has no fixed points and $\theta_i$ increases with time,
showing periodic oscillations with frequency $\sqrt{\eta_i+I_i}/\pi$~\cite{erm96}. The bifurcation
at $\eta_i+I_i=0$ is a SNIC bifurcation. Note that under the transformation $V=\tan{(\theta/2)}$
a network of theta neurons is exactly equivalent to a network of quadratic integrate-and-fire
neurons with infinite threshold and reset values~\cite{bor17}. We now proceed to analyse the network
dynamics, using ideas similar to those
in~\cite{chahat17,blasche2020,laibla20,resott14,laibla21}.

We assume that the network is characterised by two functions: firstly
the degree distribution $P(\bk)$, normalised
to sum to 1, where
$\bk=(k_{in},k_{out})$ and $k_{in}$ and $k_{out}$ are the in- and out-degrees of a neuron with
degree $\bk$, respectively, and secondly the assortativity function $a(\bk'\to\bk)$ giving the
probability of a connection from a neuron with degree $\bk'$ to one with degree $\bk$, given that
such neurons exist. We also make the mean field assumption 
that the dynamics of a neuron depend only on its degree $\bk$,
thus effectively averaging the dynamics of all neurons with the same degree.

In the limit of large $N$ and large in- and out-degrees, the network is described by the
distribution $f(\theta,\eta|\bk,t)$ where $f(\theta,\eta|\bk,t)d\theta d\eta$
is the probability that a neuron with degree $\bk$ has phase in $[\theta,\theta+d\theta]$ and
value of $\eta$ in $[\eta,\eta+d\eta]$ at time $t$. This distribution satisfies the continuity
equation
\be
   \frac{\partial f}{\partial t}+\frac{\partial}{\partial\theta}(vf)=0
\ee
where (from~\eqref{eq:dthetadt}-\eqref{eq:dudt})
\be
   v(\theta,\bk,\eta,t)=1-\cos{\theta}+(1+\cos{\theta})[\eta+I(\bk,t)],
\ee
\be
  I(\bk,t)=\frac{KN}{\mk}\sum_{\bk'}P(\bk')a(\bk'\to\bk)u(\bk',t),
\ee
and
\be
   \tau\frac{du(\bk,t)}{dt}=\widehat{F}(\bk,t)-u(\bk,t),
\ee
where $\widehat{F}(\bk,t)$ is the firing rate of neurons with degree $\bk$ at time $t$.

The Ott/Antonsen ansatz gives the dynamics for the order parameter for neurons with degree 
$\bk$~\cite{laibla20,blasche2020,lai14A}:
\begin{gather*}
   \frac{\partial b(\bk,t)}{\partial t}=\frac{-i[b(\bk,t)-1]^2}{2}+\frac{[b(\bk,t)+1]^2}{2}\\
 \times \left[-\Delta+i\eta_0+i\frac{KN}{\mk}\sum_{\bk'}P(\bk')a(\bk'\to\bk)u(\bk',t)\right] 
\end{gather*}
where $b(\bk,t)$ is the expected value of $e^{i\theta}$ for neurons with degree $\bk$, i.e.
\be
  b(\bk,t)=\int_{-\infty}^\infty\int_{0}^{2\pi}f(\theta,\eta|\bk,t)e^{i\theta}d\theta\ d\eta.
\ee
The firing rate of neurons with degree $\bk$ at time $t$ is the expected value of the
flux through $\theta=\pi$~\cite{lai15,monpaz15}, i.e.
\be
   \widehat{F}(\bk,t)=\int_{-\infty}^\infty f(\pi,\eta|\bk,t)d\eta=\frac{1}{\pi}\mbox{Re}\left(\frac{1-\bar{b}(\bk,t)}{1+\bar{b}(\bk,t)}\right)
\ee
where the overline indicates conplex conjugate. We define $F(b(\bk,t))\equiv \widehat{F}(\bk,t)$.

With neutral assortativity~\cite{resott14}, 
\be
   a(\bk'\to\bk)=\frac{k_{out}'k_{in}}{N\mk}
\ee
and with independent in- and out-degrees the degree distribution $P(\bk')$
factorises as $P(\bk')=p_{in}(k_{in}')p_{out}(k_{out}')$, where $p_{in}$ and $p_{out}$
are the marginal distributions of the relevant degrees, so
\begin{gather}
  \frac{KN}{\mk}\sum_{\bk'}P(\bk')a(\bk'\to\bk)u(\bk',t) \label{eq:sum} \\
 =\frac{K k_{in}}{\mk^2}\sum_{k_{in}'}\sum_{k_{out}'}p_{in}(k_{in}')p_{out}(k_{out}')k_{out}'u(k_{in}',k_{out}',t). \nonumber
\end{gather}
This quantity is independent of $k_{out}$ and contributes to
 the ``input'' to neurons with degree $\bk$. Thus
$b(\bk,t)$ must also be independent of $k_{out}$ and so must $F(b(\bk,t))$ and $u(\bk,t)$.
Thus~\eqref{eq:sum} simplifies to
\begin{gather}
   \frac{KN}{\mk}\sum_{\bk'}P(\bk')a(\bk'\to\bk)u(\bk',t) \nonumber \\
=\frac{K k_{in}}{\mk}\sum_{k_{in}'}p_{in}(k_{in}')u(k_{in}',t)
\end{gather}
and we see that the distribution of out-degrees does not affect the expected dynamics.
This was observed by Roxin in~\cite{rox11}, although he observed that broadening the out-degree
distribution increases the amplitude of the cross-correlation of synaptic currents, something
we do not consider here.

We have 
\be
   \tau\frac{du(k_{in},t)}{dt}=F(b(k_{in},t))-u(k_{in},t),
\ee
and defining 
\be
   s(t)\equiv\sum_{k_{in}}p_{in}(k_{in})u(k_{in},t)
\ee
we see that $s$ satisfies
\be
   \tau\frac{ds}{dt}=\sum_{k_{in}}p_{in}(k_{in})F(b(k_{in},t))-s \label{eq:dsdt}
\ee
and the dynamics of $b$ are given by
\begin{gather}
   \frac{\partial b(k_{in},t)}{\partial t}=\frac{-i[b(k_{in},t)-1]^2}{2} \nonumber \\
+\frac{[b(k_{in},t)+1]^2}{2} \left[-\Delta+i\eta_0+\frac{i K k_{in} s}{\mk}\right]. \label{eq:dbdt}
\end{gather}
Equations~\eqref{eq:dsdt}-\eqref{eq:dbdt} form a set of $N_{k_{in}}+1$ ordinary differential
equations (ODEs) governing the network's dynamics, where $N_{k_{in}}$ is the number of distinct
in-degrees in the network. In the next section we give some numerical results showing the
possible dynamics of this set of equations.

\subsection{Results}
We first consider inhibitory coupling, i.e.,~$K<0$.

\subsubsection{Inhibitory coupling}
\begin{figure}
\begin{center}
\includegraphics[width=3.0in]{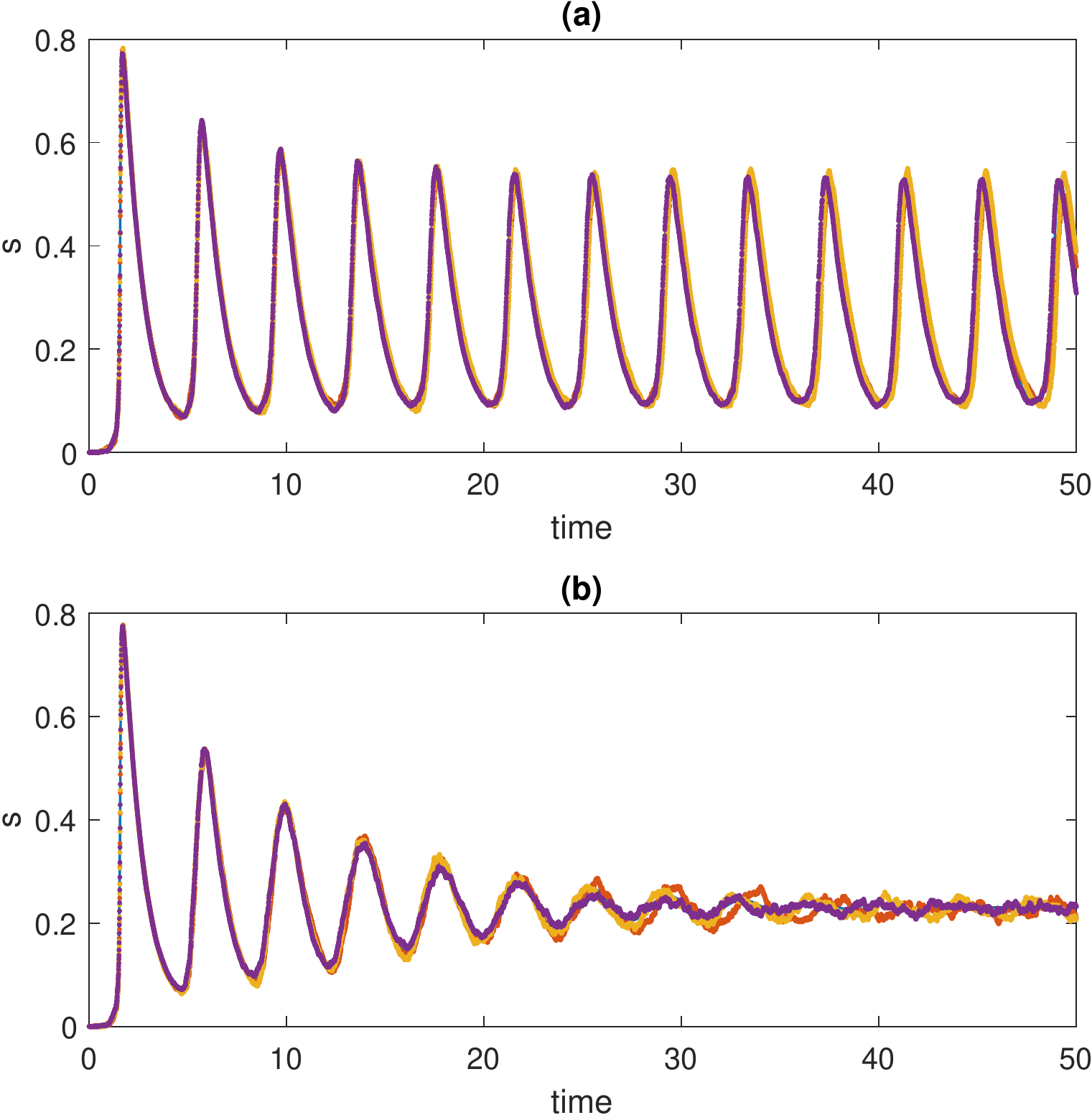}
\caption{(a) $s$ for the reduced model~\eqref{eq:dsdt}-\eqref{eq:dbdt} (blue curve) and
for the full model~\eqref{eq:dthetadt}-\eqref{eq:dudt} (dots) for $p_{out}$ being uniform
on $[10,190],[50,150]$ and $[90,110]$ (different colours). $p_{in}$ is uniform on $[95,105]$
(i.e.~$\sigma=5$).
A different realisation of the
$\eta_i$ was used for the different networks. (b): as for (a) but 
now $p_{in}$ is uniform on $[50,150]$
($\sigma=50$). For~~\eqref{eq:dthetadt}-\eqref{eq:dudt}
$N=500$ neurons were used and the initial conditions were $\theta_i=u_i=0$,
and for~\eqref{eq:dsdt}-\eqref{eq:dbdt} we used $b=1$ and $s=0$.}
\label{fig:exam}
\end{center}
\end{figure}

Consider the parameter values $\eta_0=1,\Delta=0.05,\tau=1,K=-2,\mk=100$. Having $\eta_0>0$
indicates that when uncoupled, most neurons would be firing rather than quiescent.
We choose $p_{in}$
to be uniform with mean $\mk=100$ and write its support as $[100-\sigma,100+\sigma]$. 
For $\sigma=5$ (a narrow distribution) we obtain global oscillations, see Fig.~\ref{fig:exam}(a).
However, when the in-degree distribution is made broader ($\sigma=50$), the oscillations die out: 
see Fig.~\ref{fig:exam}(b). Note the independence of the dynamics on the distribution of
out-degrees, $p_{out}(k_{out})$, as expected.

To numerically solve~\eqref{eq:dsdt}-\eqref{eq:dbdt} we treat $k_{in}$ as a continuous
variable and discretise the support of $p_{in}$ using 100 evenly spaced points, and use
$p_{in}=1/100$ at each of those points, effectively using the midpoint rule. 
To create the network used in~\eqref{eq:dthetadt}-\eqref{eq:dudt} we
randomly sample $N$ in-degrees from $p_{in}(k_{in})$ and $N$ out-degrees from $p_{out}(k_{out})$, 
choosing until the sum of the in-degrees
equals the sum of the out-degrees, then use the configuration model 
to connect the network~\cite{new03}.
Any self or multiple connections are then removed by random rewiring, keeping the degrees
fixed. To solve~\eqref{eq:dthetadt}-\eqref{eq:dudt} we used Euler's method with a stepsize
of $0.001$.

\begin{figure}
\begin{center}
\includegraphics[width=3.0in]{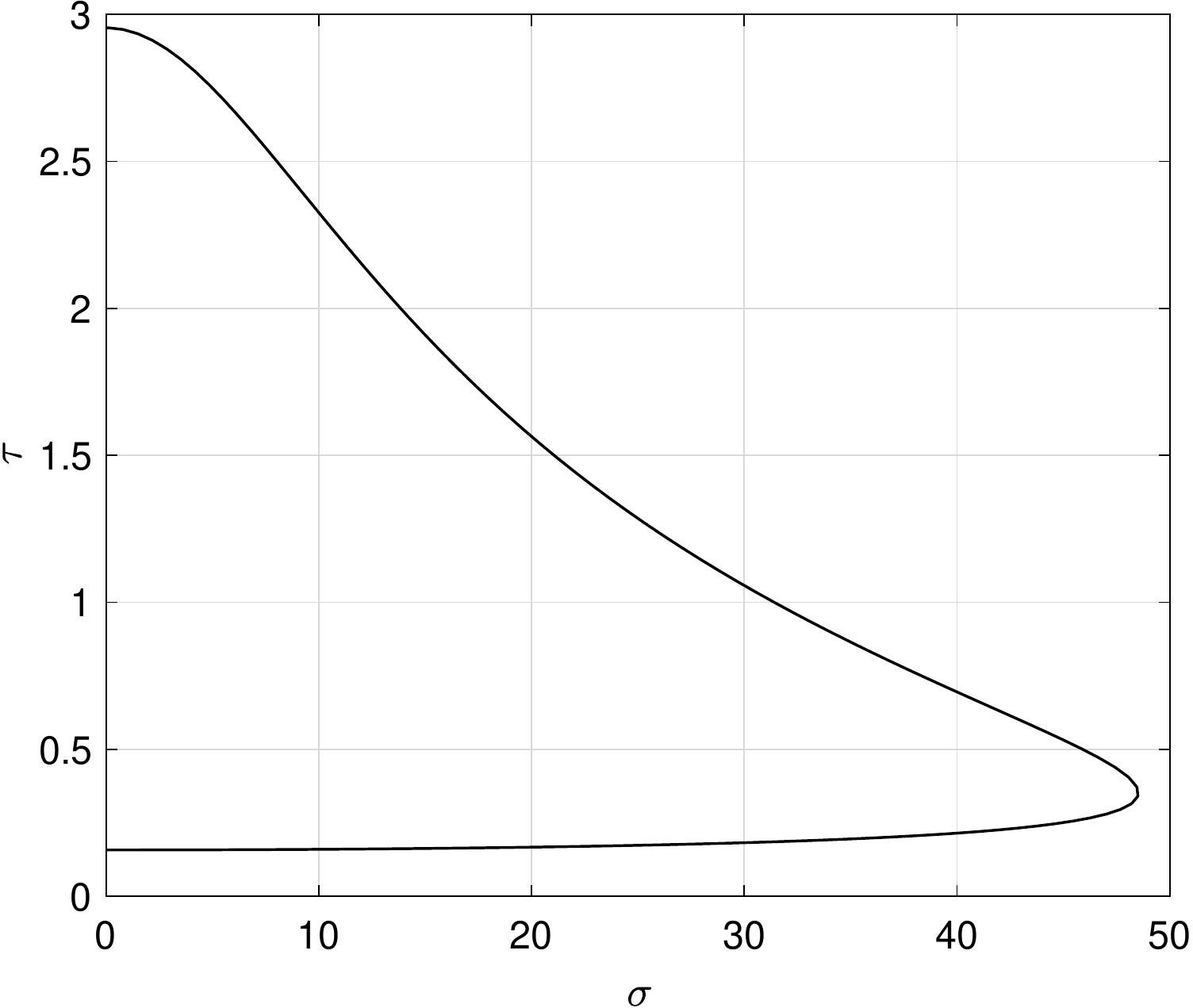}
\caption{ 
Hopf bifurcation curve for a fixed point of~\eqref{eq:dsdt}-\eqref{eq:dbdt}. 
A stable periodic orbit exists to the left of the curve and a stable fixed point to the right.
Other parameters: $\eta_0=1,\Delta=0.05,K=-2,\mk=100$. $p_{in}$ is uniform on 
$[100-\sigma,100+\sigma]$.}
\label{fig:sigtau}
\end{center}
\end{figure}

The destruction of oscillations seen in Fig.~\ref{fig:exam}
seems due to a Hopf bifurcation. Using pseudo-arclength continuation to follow the stable
fixed point of~\eqref{eq:dsdt}-\eqref{eq:dbdt} as $\sigma$ is decreased we find a Hopf
bifurcation at $\sigma\approx 31.4$, and continuing that bifurcation as both $\sigma$ and
$\tau$ are varied we obtain the curve in Fig.~\ref{fig:sigtau}. For any $\tau$ for which
an oscillation occurs, increasing $\sigma$ will destroy the oscillations, and for small $\sigma$,
a value of $\tau$ which is either too large or too small will also destroy oscillations.
(We varied $\tau$ here just as an example; we could equally well 
vary other parameters such as $\eta_0$ or
$\Delta$.)
The destruction of oscillations in an inhibitory network by broadening the in-degree distribution
was also observed by Roxin~\cite{rox11}. He analysed a heuristic rate model containing a fixed delay
(since he used delayed synapses) and found a Hopf bifurcation in that model, in agreement with
the results shown here.

\begin{figure}
\begin{center}
\includegraphics[width=3.0in]{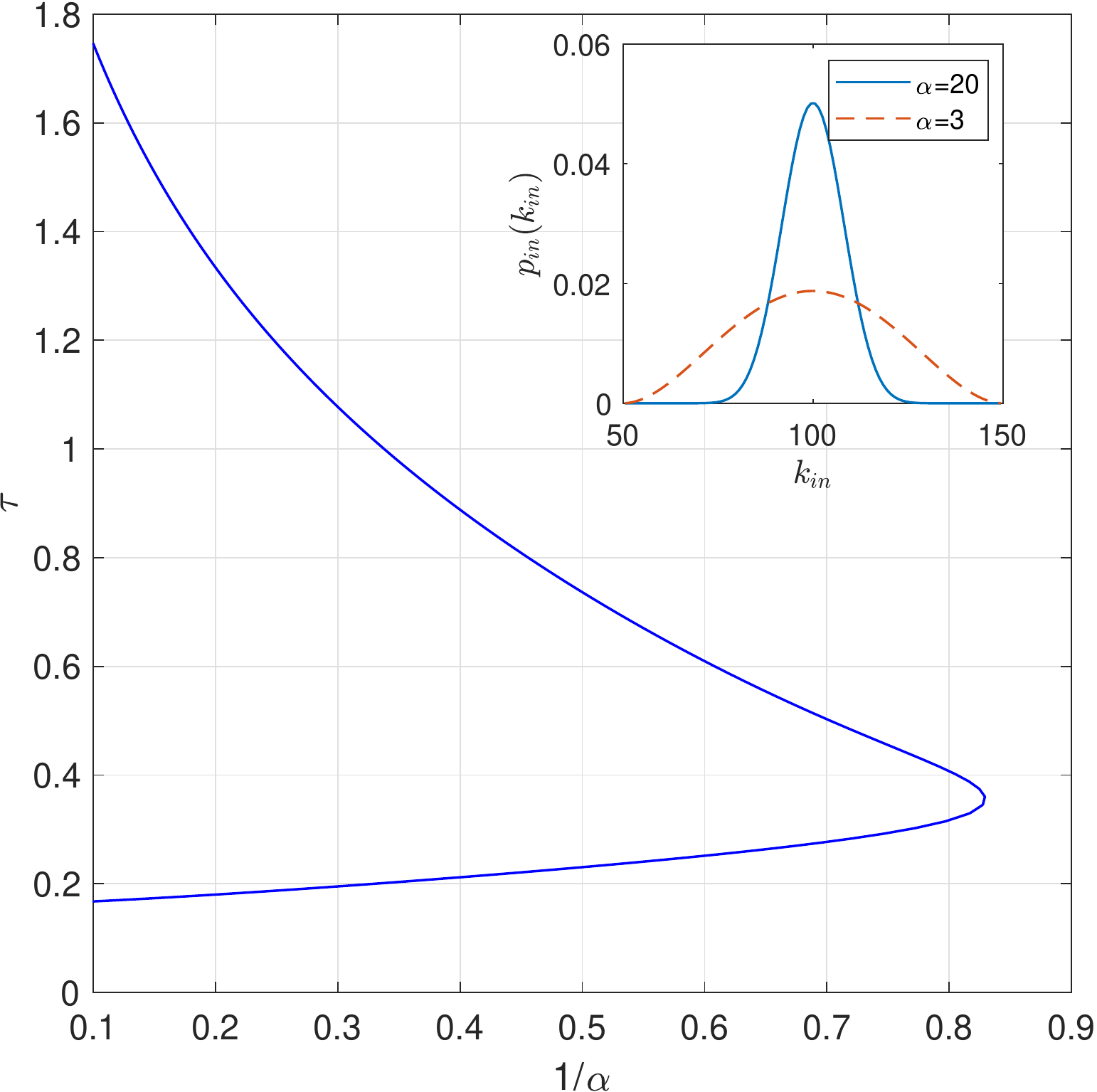}
\caption{A beta distribution of in-degrees. 
Hopf bifurcation curve for fixed point of~\eqref{eq:dsdt}-\eqref{eq:dbdt}. 
A stable periodic orbit exists to the left of the curve and a stable fixed point to the right.
Other parameters: $\eta_0=1,\Delta=0.05,K=-2,\mk=100$. Inset shows the beta distribution
on $[50,150]$ with $\alpha=3,20$.}
\label{fig:alptau}
\end{center}
\end{figure}

To investigate the generality of our result
we now consider a beta distribution of in-degrees,
with equal parameters greater than one, shifted to
have mean 100 and support on $[50,150]$, i.e.
\be
   p_{in}(k_{in})=\begin{cases} Cx^{\alpha-1}(1-x)^{\alpha-1}, & 0\leq x \leq 1 \\ 0, & \mbox{otherwise} \end{cases}
\ee
where $x=(k_{in}-50)/100$ and $C$ is a normalisation factor. Increasing $\alpha$ narrows the
distribution, as shown in the inset of Fig.~\ref{fig:alptau}. Varying $\alpha$ and $\tau$
we find a curve of Hopf bifurcations, shown in Fig.~\ref{fig:alptau}, which shows the same
qualitative behaviour as for the uniform distribution. 


\begin{figure}
\begin{center}
\includegraphics[width=3.0in]{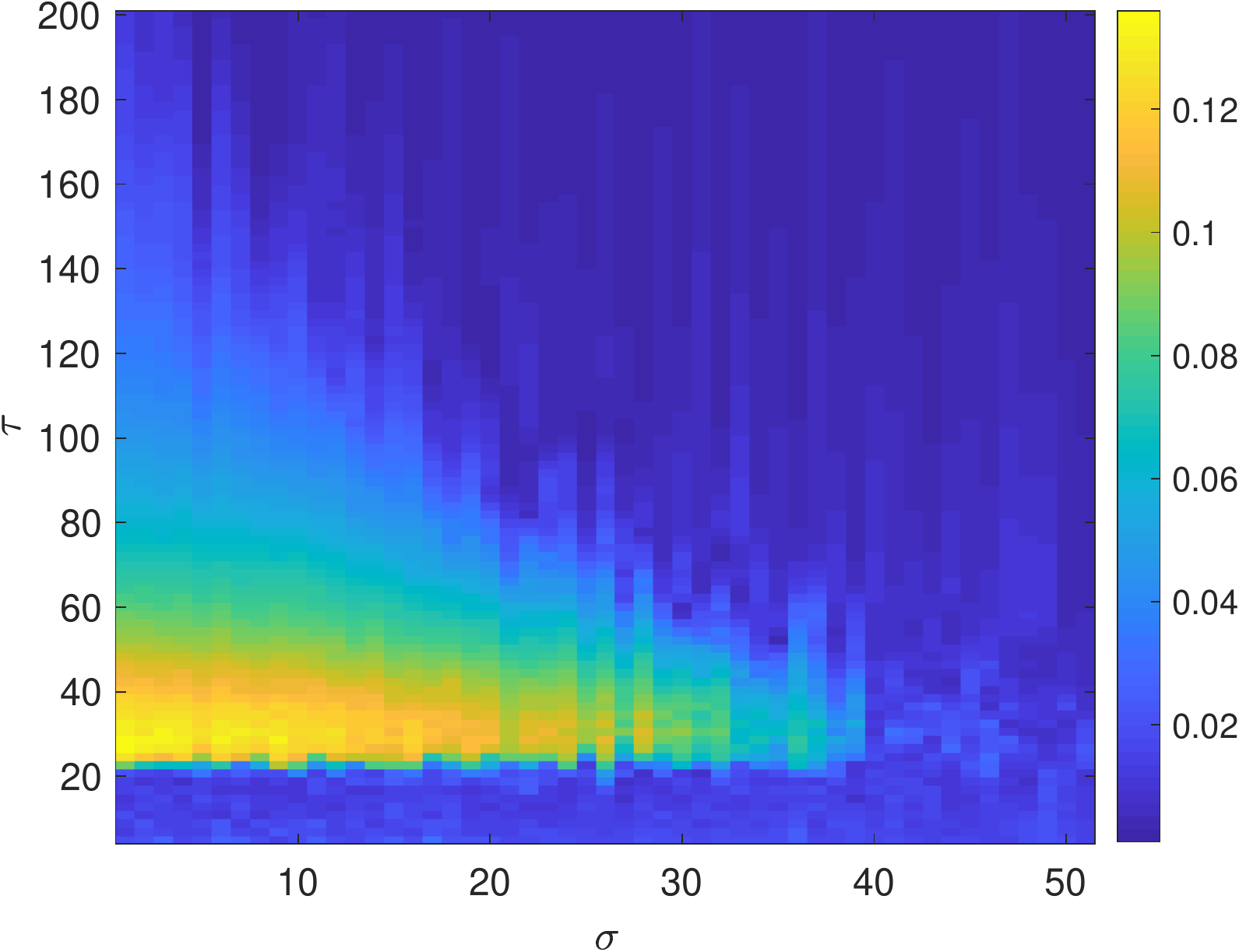}
\caption{Results for a Morris-Lecar model. $\tau$ is the synaptic timescale and the
in-degree distribution is uniform on $[100-\sigma,100+\sigma]$.
Colour shows the standard deviation of $\hat{s}$ over 5 seconds of simulated time, having
discarded the first 45 seconds. The same $\eta_i$ were used for each simulation.
We used $N=500$ neurons.}
\label{fig:ML}
\end{center}
\end{figure}

To further demonstrate the generality of our results we
now consider a network of 500 Morris-Lecar neurons~\cite{tsukit06}, 
known to undergo a SNIC bifurcation as the
input current is increased~\cite{laibla20}.
The network equations are given in Appendix~\ref{sec:app}.
The in-degree distribution is uniform on $[100-\sigma,100+\sigma]$ and the 
out-degree is uniform on $[50,150]$. For each different value of $\sigma$ we generate a
network as explained above and for each network we vary $\tau$, the synaptic timescale,
integrating for $50$ seconds
at each value of $\tau$.
Defining $\hat{s}=N^{-1}\sum_{i=1}^N s_i$
we discard data from the first $45$ seconds and calculate the standard deviation of $\hat{s}$
over the last 5 seconds, plotting that in Fig.~\ref{fig:ML}. Large values indicate
oscillations while small values indicate an approximate steady state. We see results
consistent with those in Figs.~\ref{fig:sigtau} and~\ref{fig:alptau}.
Thus we conclude that broadening the in-degree distribution of an inhibitory network
of Type-I neurons acts to destroy global oscillatory behaviour. This is presumably due to
having a wider range of dynamics for neurons with different in-degrees, making them harder
to synchronise.
We now consider excitatory coupling, i.e.,~$K>0$.

%
%
%
%

\subsubsection{Excitatory coupling}
\begin{figure}
\begin{center}
\includegraphics[width=3.0in]{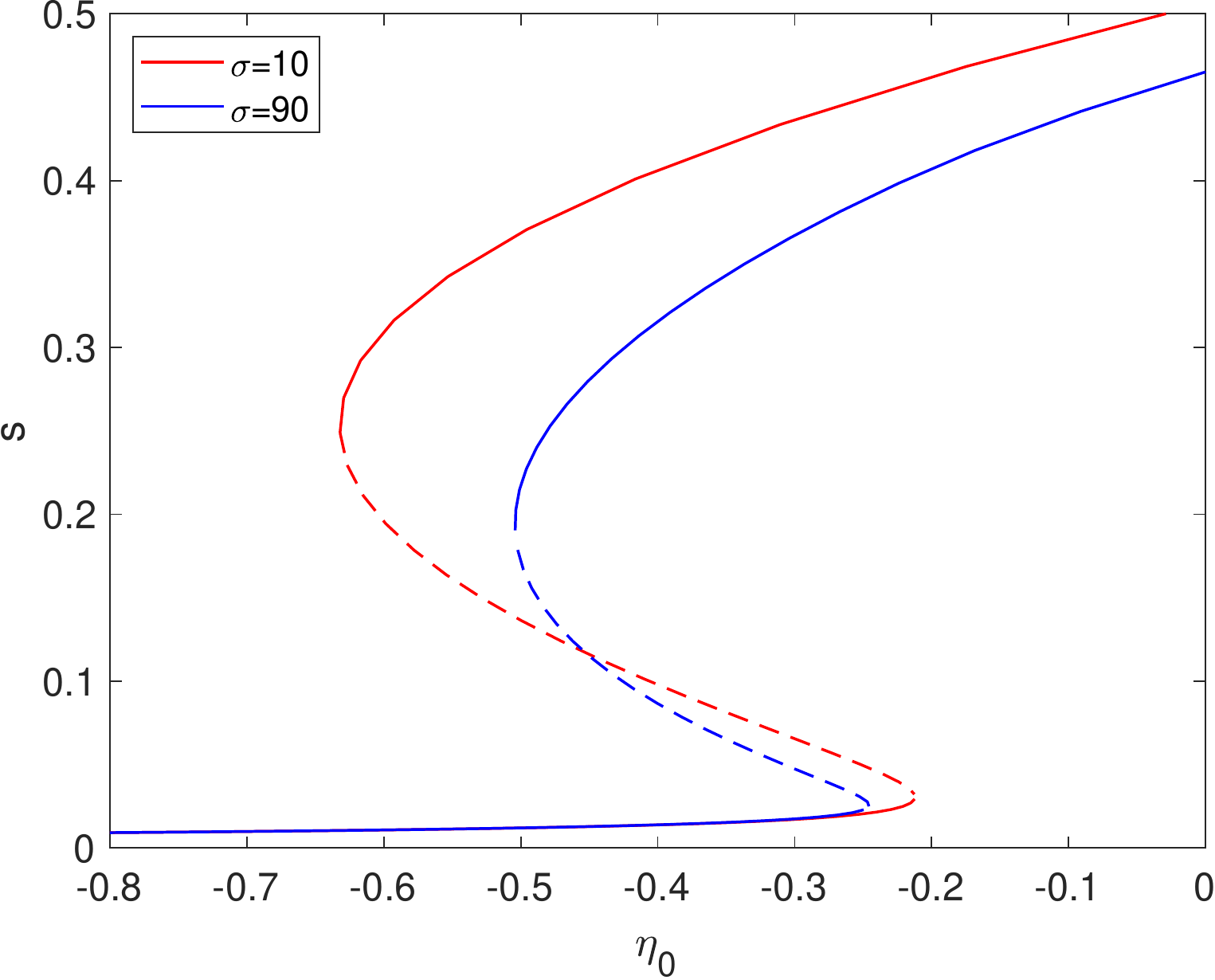}
\caption{$s$ at fixed points of~\eqref{eq:dsdt}-\eqref{eq:dbdt} for $\sigma=10$ (red)
and $\sigma=90$ (blue). Solid curves are stable, dashed unstable.
Parameters: $\Delta=0.05,\tau=1,K=5,\mk=100$, uniform in-degree on $[100-\sigma,100+\sigma]$.}
\label{fig:excit}
\end{center}
\end{figure}

\begin{figure}
\begin{center}
\includegraphics[width=3.0in]{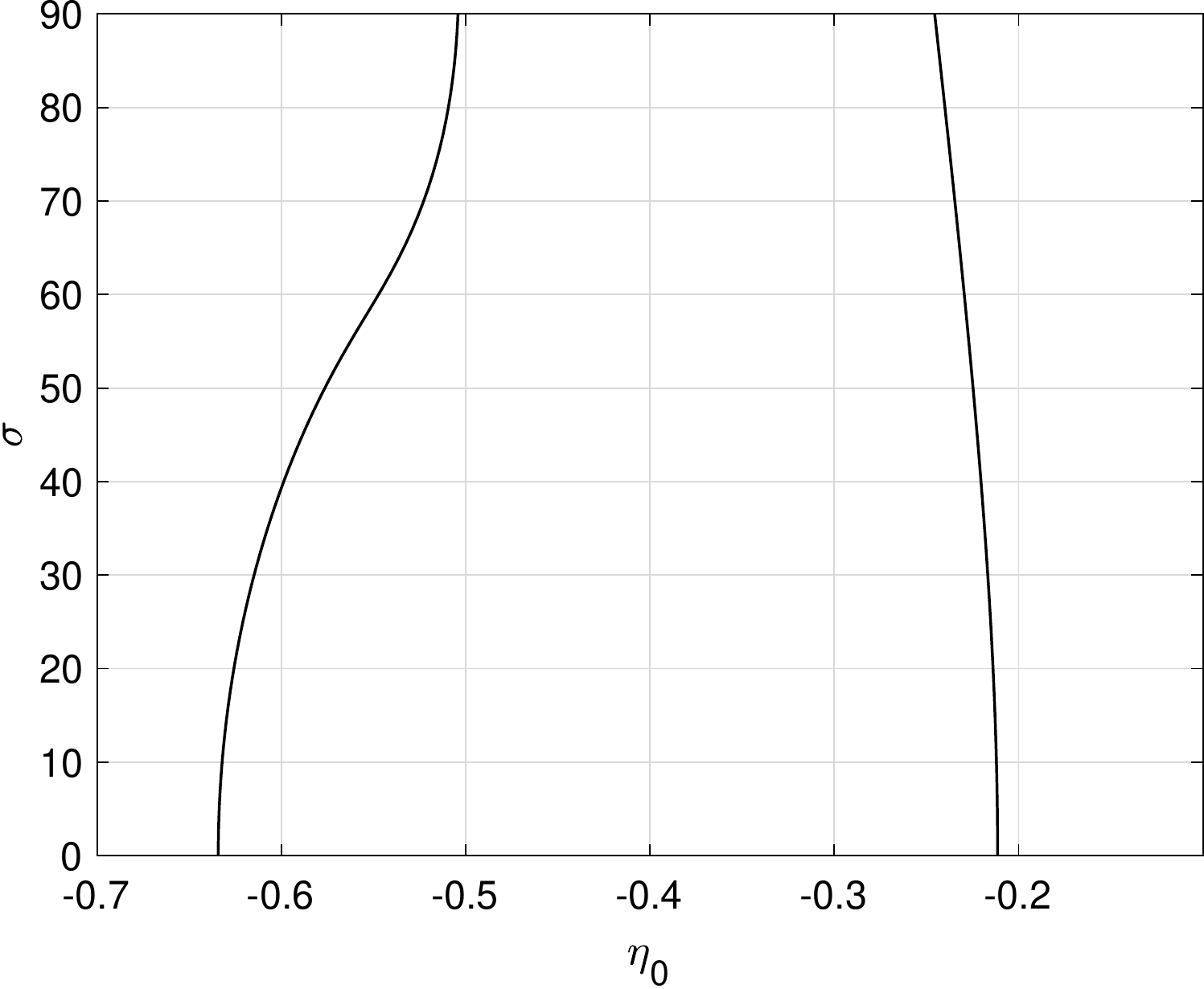}
\caption{Curves of saddle-node bifurcations of fixed points of~\eqref{eq:dsdt}-\eqref{eq:dbdt}.
The network is bistable between the curves and has a single attractor outside this region.
Fig.~\ref{fig:excit} corresponds to horizontal ``slices'' through this figure at $\sigma=10$ and
$\sigma=90$.
Parameters: $\Delta=0.05,\tau=1,K=5,\mk=100$, uniform in-degree on $[100-\sigma,100+\sigma]$.}
\label{fig:sn}
\end{center}
\end{figure}

Consider the parameter values $\Delta=0.05,\tau=1,K=5,\mk=100$. Varying $\eta_0$ we expect
a region of bistability between a high activity steady state and a low activity steady state,
as is often found in excitatory networks~\cite{lailon03}.
This is found, as shown in Fig.~\ref{fig:excit}, where saddle-node bifurcations
mark the boundaries of the bistable region. Varying the
width of the in-degree distribution ($\sigma$) varies the width of the bistable region.
In particular, widening the in-degree distribution narrows the width of the bistable region.
Following the saddle-node bifurcations seen in Fig.~\ref{fig:excit} we obtain
Fig.~\ref{fig:sn}. Similar behaviour was seen for a beta distribution of in-degrees (not shown).

\begin{figure}
\begin{center}
\includegraphics[width=3.0in]{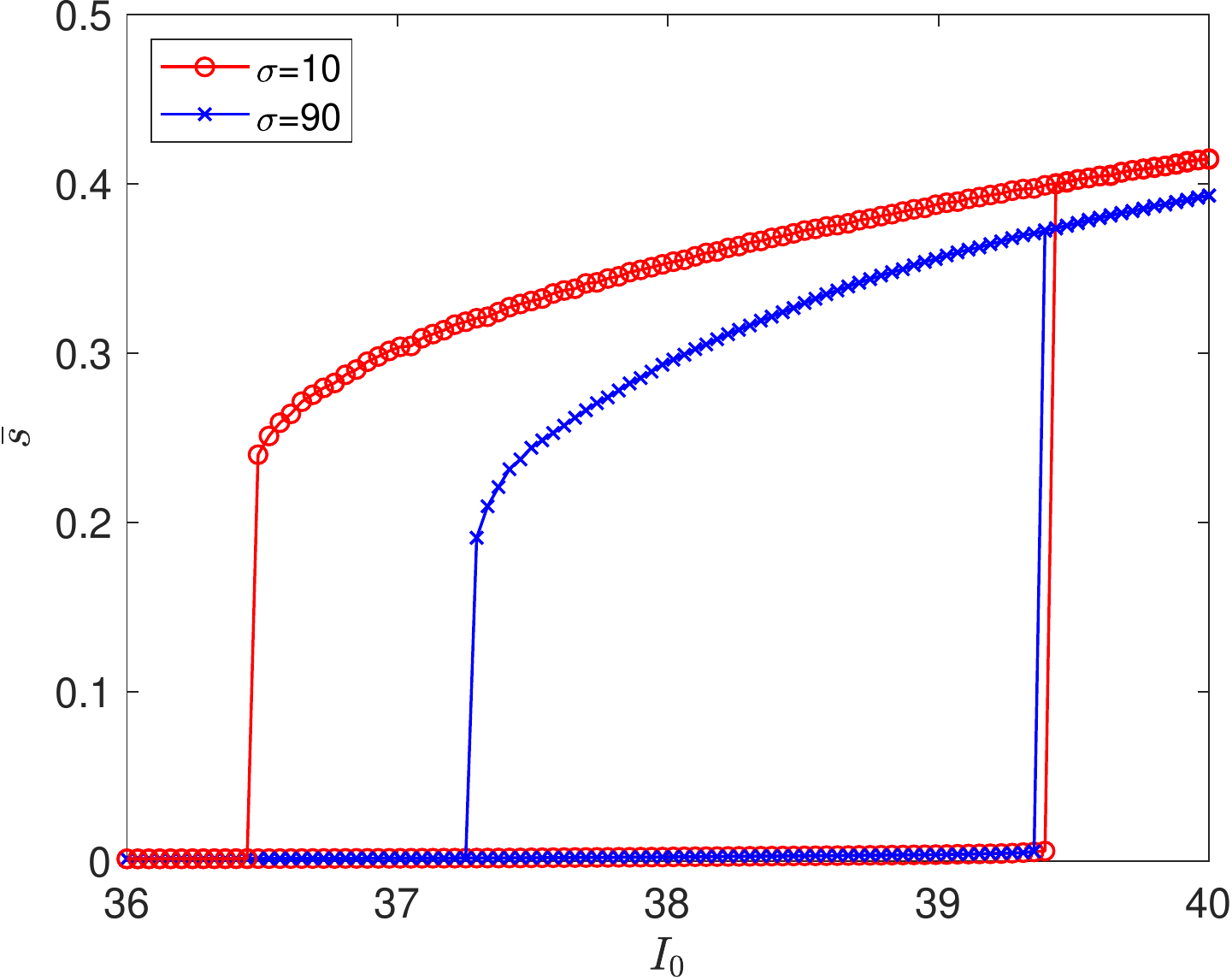}
\caption{Approximate steady states for a Morris-Lecar model.
For each of the two networks $I_0$ was quasistatically increased up and then down.
Both networks show bistability for a range of $I_0$ values, and vertical jumps at
apparent saddle-node bifurcations. Compare with Fig.~\ref{fig:excit}.
Parameters: $\mk=100$, uniform in-degree on $[100-\sigma,100+\sigma]$, out-degree is uniform
on $[50,150]$. The same $\eta_i$ were used for each simulation.}
\label{fig:MLex}
\end{center}
\end{figure}

We reproduced this behaviour in a network of $N=500$ Morris-Lecar neurons with parameters
as given in  Appendix~\ref{sec:app}.  For networks with $\sigma=10$ and $90$ 
we quasistatically varied $I_0$, integrating
for $10$ seconds at each value. We define $\bar{s}$ to be the mean of $\hat{s}$ over the
last $2$ seconds of simulation and plot this in Fig.~\ref{fig:MLex}. The results are qualitatively
the same as in Fig.~\ref{fig:excit}: $\bar{s}$ is lower when $\sigma=90$ than when $\sigma=10$,
and the left-most saddle-node bifurcation is moved more than the right-most when $\sigma$ is
varied. The threshold for firing for single neuron is $I_0\approx 39.69$ so the jumps
occur at values of $I_0$ less than this, consistent with the results in Fig.~\ref{fig:sn}.
Thus we conclude that for an excitatory network of Type-I neurons, broadening the in-degree
distribution narrows the range of values of the mean input for which the network
is bistable. The effects of varying other parameters could equally well be investigated 
using the techniques shown here. 

\section{Gap junctions}
\label{sec:gap}

We now consider theta neurons coupled by gap junctions~\cite{Benzuk04,galhes99}. 
Gap junctional coupling is well-known to induce synchrony in networks
of neurons\cite{ostbru09,erm06,trakop01}. The 
quadratic integrate-and-fire (QIF) neuron~\cite{latric00} 
with infinite firing threshold and reset to
$V=-\infty$ is equivalent under the transformation $V=\tan{(\theta/2)}$ to a theta 
neuron~\cite{ermkop86}, and
since gap junction coupling is through voltage differences it is easier to start with a 
network of QIF neurons. Our analysis is similar to that in~\cite{lai15}; also 
see~\cite{monpaz20,piedev19,Byrros20}.
A theta neuron is an excitable system, so our results add to those on coupled excitable
systems~\cite{lafcol10,shikur86,chistr08,zhepik19}.

\subsection{Model and Theory}
Consider a network of $N$ gap-junction coupled QIF neurons
governed by
\be
  \frac{dV_j}{dt}=\eta_j+V_j^2+\frac{g}{\mk}\sum_{l=1}^N A_{jl}(V_l-V_j) \label{eq:dVdta}
\ee
for $j=1,\dots N$ together with the rule that if $V_j(t^-)=\infty$ then $V_j(t^+)=-\infty$, and neuron
$j$ is said to fire at this time $t$. $A$ describes the connectivity of the network, 
where $A_{jl}=1$ if neurons $j$ and $l$ are connected and zero otherwise. 
Since gap junctional coupling is not directional we have $A_{jl}=A_{lj}$.
$k_j$ is the degree of the $j$th neuron, i.e.~$k_j=\sum_l A_{jl}$, $\mk$ is the mean degree,
as above, and $g$ is the strength of coupling (non-negative). The $\eta_j$ are 
randomly chosen from
a distribution $h(\eta)$.

We 
rewrite~\eqref{eq:dVdta} as
\be
  \frac{dV_j}{dt}=\eta_j+V_j^2-\frac{gk_jV_j}{\mk}+\frac{g}{\mk}\sum_{l=1}^N A_{jl}V_l.
\ee
Now let $V_j=\tan{(\theta_j/2)}$. Then
\begin{align}
  \frac{dV_j}{dt}  & =  \frac{d\theta_j/dt}{2\cos^2{(\theta_j/2)}} \nonumber \\
 & =  \eta_j+\tan^2{(\theta_j/2)}-\frac{gk_j}{\mk}\tan{(\theta_j/2)} \nonumber \\
 &  +\frac{g}{\mk}\sum_{l=1}^N A_{jl}\tan{(\theta_l/2)}
\end{align}
so
\begin{gather}
   \frac{d\theta_j}{dt}=1-\cos{\theta_j}+(1+\cos{\theta_j}) \nonumber \\
\times \left[\eta_j-\frac{gk_j}{\mk}\tan{(\theta_j/2)}+\frac{g}{\mk}\sum_{l=1}^N A_{jl}\tan{(\theta_l/2)}\right].
\label{eq:dth1}
\end{gather}
Noting that
\be
   \tan{(\theta/2)}=\frac{\sin{\theta}}{1+\cos{\theta}},
\ee
we have
\begin{gather}
   \frac{d\theta_j}{dt}=1-\cos{\theta_j}-\frac{gk_j}{\mk}\sin{\theta_j}+(1+\cos{\theta_j}) \nonumber \\
\times \left[\eta_j+\frac{g}{\mk}\sum_{l=1}^N A_{jl}\tan{(\theta_l/2)}\right].
\label{eq:dth2}
\end{gather}
When a neuron fires, at $\theta=\pi$, the term involving $\tan$ becomes infinite.
To avoid this problem we follow~\cite{erm06} 
and replace $\tan{(\theta/2)}$ in~\eqref{eq:dth2} by 
\be
   q(\theta)\equiv \frac{\sin{\theta}}{1+\cos{\theta}+\epsilon},
\ee
where $0<\epsilon\ll 1$, thereby removing the singularity. We take the limit $\epsilon\to 0$ below.

We analyse
the system in a similar way as in Sec.~\ref{sec:syn}.
The system is described by the probability density function $f(\eta,\theta|{\bf k},t)$ which satisfies~\cite{ome14,str00,abrmir08}
\be
  \frac{\partial f}{\partial t}+\frac{\partial}{\partial \theta}(fv)=0, \label{eq:cont}
\ee
where 
\begin{gather}
   v(\eta,\theta,{\bf k},t)\equiv 1-\cos{\theta}-\frac{gk_{in}}{\mk}\sin{\theta} \nonumber \\
+(1+\cos{\theta})\left[\eta+gT({\bf k},t)\right], \label{eq:v}
\end{gather}
where
\be
   T({\bf k},t)\equiv \frac{N}{\mk}\sum_{{\bf k}'}P({\bf k}')a({\bf k}'\rightarrow{\bf k})Q({\bf k}',t),
\ee
\be
   Q({\bf k}',t)=\int_{-\infty}^\infty\int_0^{2\pi}f(\eta,\theta|{\bf k}',t)q(\theta)\ d\theta \ d\eta
\ee
and ${\bf k}=(k_{in},k_{out})$, $P({\bf k}')$ is the distribution of degrees of a neuron, and
$a({\bf k}'\rightarrow{\bf k})$ is the probability of a 
connection from a neuron with degree ${\bf k}'$
to one with degree ${\bf k}$. But connections are undirected, so a neuron just has a degree.
Thus we write
\be
   T(k,t)= \frac{N}{\mk}\sum_{k'}P(k')a(k'\rightarrow k)Q(k',t)
\ee
where $P(k)$ is the degree distribution, and
\be
   Q(k',t)=\int_{-\infty}^\infty\int_0^{2\pi}f(\eta,\theta|k',t)q(\theta)\ d\theta \ d\eta
\ee
is the expected value of $q(\theta)$ for neurons with degree $k'$ at time $t$.
From~\cite{lai15}, assuming that $h(\eta)$ is a Lorentzian with median $\eta_0$
and HWHM $\Delta$, we have
\begin{gather}
   \frac{\partial b(k,t)}{\partial t}=\frac{(i\eta_0-\Delta)(1+b)^2-i(1-b)^2}{2} \nonumber \\
 +\frac{i(1+b)^2gT+g(1-b^2)}{2} \label{eq:dzdt}
\end{gather}
where $T(k,t)$ and $Q(k',t)$ are as above and
\be
   b(k,t)=\int_{-\infty}^\infty\int_0^{2\pi}f(\eta,\theta|k,t)e^{i\theta}\ d\theta \ d\eta
\ee
is the complex-valued order parameter for neurons with degree $k$ at time $t$.

By expanding $q(\theta)$ in a Fourier series it was shown in~\cite{lai15} that  
\be 
  Q(k',t)=\sum_{m=1}^\infty \left[c_m b(k',t)^m+\mbox{c.c.}\right] \label{eq:Q}
\ee
where ``c.c.'' is the complex conjugate of the previous term and
\be
   c_m=\frac{i[(\sqrt{2\epsilon+\epsilon^2}-1-\epsilon)^{m+1}-(\sqrt{2\epsilon+\epsilon^2}-1-\epsilon)^{m-1}]}{2\sqrt{2\epsilon+\epsilon^2}}.
\ee
Recently,~\cite{piedev19} showed that
\be
   \lim_{\epsilon\rightarrow 0} c_m=i(-1)^m
\ee
and thus
\begin{gather}
   Q(k',t)=i\sum_{m=1}^\infty (-1)^m[b(k',t)^m-\bar{b}(k',t)^m] \nonumber \\
 =\frac{2\mbox{Im}[b(k',t)]}{[1+b(k',t)][1+\bar{b}(k',t)]}
\end{gather}
where we have summed the geometric series.

Defining
\be
   w\equiv\frac{1-\bar{b}}{1+\bar{b}} \label{eq:w}
\ee
(so $b=(1-\bar{w})/(1+\bar{w})$) we find that $w$ satisfies 
\be
   \frac{\partial w(k,t)}{\partial t}=i\eta_0+\Delta-i[w(k,t)]^2+igT(k,t)-gw(k,t), \label{eq:dwdt}
\ee
and writing $w=\pi \phi+iV$ where $\phi$ and $V$ are real
we find that $Q(k',t)=V(k',t)$ and the real and imaginary parts of~\eqref{eq:dwdt} give
\begin{align}
   \frac{\partial \phi(k,t)}{dt} & =\frac{\Delta}{\pi}+2\phi(k,t)V(k,t)-g\phi(k,t) \label{eq:dfdt} \\
   \frac{\partial V(k,t)}{dt} & =\eta_0-\pi^2[\phi(k,t)]^2+[V(k,t)]^2 \nonumber \\
 & +g[T(k,t)-V(k,t)] \label{eq:dvdt}
\end{align}
where 
\be
   T(k,t)= \frac{N}{\mk}\sum_{k'}P(k')a(k'\rightarrow k)V(k',t).
\ee
The interpretation of $\phi$ and $V$ is that $\phi(k,t)$ is the expected firing frequency
of neurons with degree $k$ at time $t$, and $V(k,t)$ is the mean voltage of QIF
neurons with degree $k$ at time $t$ where voltage $V$ and $\theta$ are related through
$V=\tan{(\theta/2)}$~\cite{monpaz15}. Note that~\eqref{eq:dfdt}-\eqref{eq:dvdt} are completely
equivalent to~\eqref{eq:dzdt}.

Assuming neutral assortativity we have 
\be
   a(k'\rightarrow k)=\frac{k'k}{N\mk}
\ee
so that
\be 
   T(k,t)= \frac{k}{\mk^2}\sum_{k'}k'P(k')V(k',t). \label{eq:T}
\ee
Note that if all neurons have the same degree then $T=V$, so the last term 
in~\eqref{eq:dvdt} vanishes, and~\eqref{eq:dfdt}-\eqref{eq:dvdt} reduce to a pair
of ODEs. A special case of this is all-to-all coupling, which was studied in~\cite{piedev19}.
Also,~\eqref{eq:T} is invariant under the scaling
$k\to\alpha k, P(k)\to P(k/\alpha)$, i.e.~only degree relative to mean degree is of
relevance.

\subsection{Results}

\begin{figure}
\begin{center}
\includegraphics[width=3.0in]{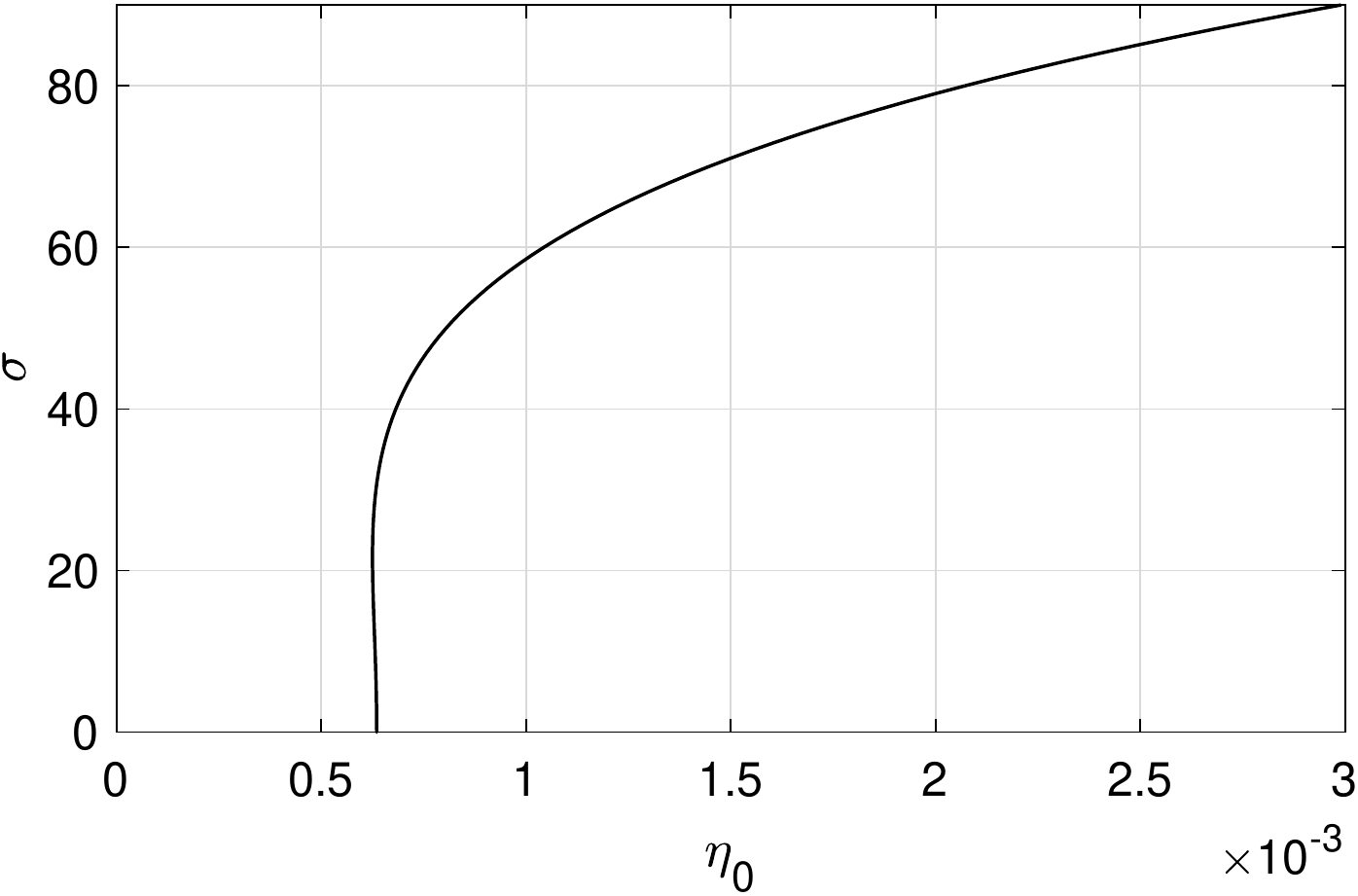}
\caption{SNIC bifurcation curve of a fixed point of~\eqref{eq:dfdt}-\eqref{eq:dvdt}. 
A stable fixed point exists to the left of the curve and stable periodic oscillations
to the right. Parameters: $\Delta=0.01,g=0.4,\mk=100$,
$P(k)$ is a uniform distribution on $[100-\sigma,100+\sigma]$.}
\label{fig:gaplargeg}
\end{center}
\end{figure}

First consider a network with $\Delta=0.01$, and $g=0.4$,
with $\mk=100$. As above we consider 
a uniform distribution of degrees  on $[100-\sigma,100+\sigma]$.
Consistent with~\cite{piedev19,monpaz20} we find that when 
increasing $\eta_0$ the transition to periodic firing is through a SNIC bifurcation, as shown
in Fig.~\ref{fig:gaplargeg}. Increasing the width of the degree distribution increases
the value of $\eta_0$ at which collective oscillations start. Note that even for $\sigma=0$
(i.e.~identical degrees) $\eta_0$ can be small and positive yet the network is quiescent,
as also found by~\cite{piedev19,monpaz20}.
Note also the small range of $\eta_0$ values as $\sigma$ is varied.

\begin{figure}
\begin{center}
\includegraphics[width=3.0in]{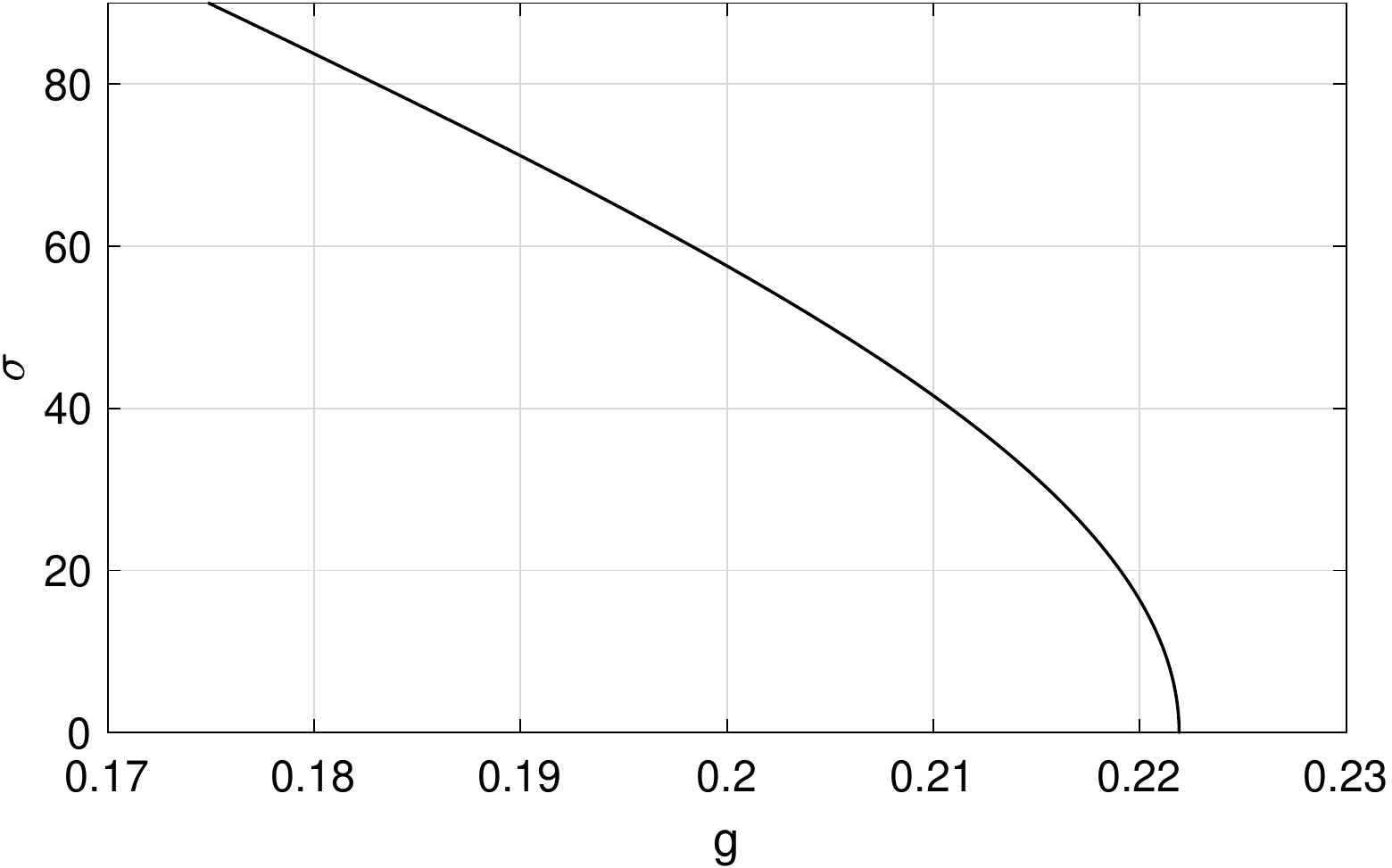}
\caption{Hopf bifurcation curve of a fixed point of~\eqref{eq:dfdt}-\eqref{eq:dvdt}. 
A stable fixed point exists to the left and stable periodic oscillations
to the right. Parameters: $\Delta=0.05,\mk=100,\eta_0=0.2$,
$P(k)$ is uniform distribution on $[100-\sigma,100+\sigma]$.}
\label{fig:gapsmallg}
\end{center}
\end{figure}

Now consider a more heterogeneous network with $\Delta=0.05$ and $\eta_0=0.2$,
i.e.~well above threshold so that most neurons would fire if uncoupled,
again with $\mk=100$ and
a uniform degree distribution on $[100-\sigma,100+\sigma]$. Consistent with the results
in~\cite{piedev19,monpaz20} we find that upon
increasing $g$ (the strength of coupling) 
the transition to firing is through a Hopf bifurcation, as shown
in Fig.~\ref{fig:gapsmallg}. Increasing the width of the degree distribution  decreases
the value of $g$ at which collective oscillations start.
This bifurcation is reminiscent of that which occurs in all-to-all
connected networks of Winfree oscillators~\cite{win67}: 
increasing the coupling strength causes the onset of oscillations through a Hopf 
bifurcation~\cite{pazmon13,galmon17,laibla21}.
As is also seen in networks of Winfree oscillators, decreasing $\Delta$ (the level of heterogeneity)
has the same effect as increasing $g$, producing oscillations via a Hopf bifurcation (not shown).


We tried to reproduce the trend in Fig.~\ref{fig:gaplargeg} using 
a network of gap junction coupled Morris-Lecar neurons, the equations of which
are given in Appendix~\ref{sec:app}. We chose the $I_i$ from a Lorentzian
with HWHM $\Delta=0.01$ and coupling strength $\epsilon=0.3$.
We were unable to reproduce the trend. This is likely due to the sensitivity of the network to the
value of $I_0$: notice the very small range of $\eta_0$ values in Fig.~\ref{fig:gaplargeg}.
The variation in values of $I_0$ at which the bifurcation occured between different networks
(with different $\sigma$) was too large to determine any significant trend.

\begin{figure}
\begin{center}
\includegraphics[width=3.0in]{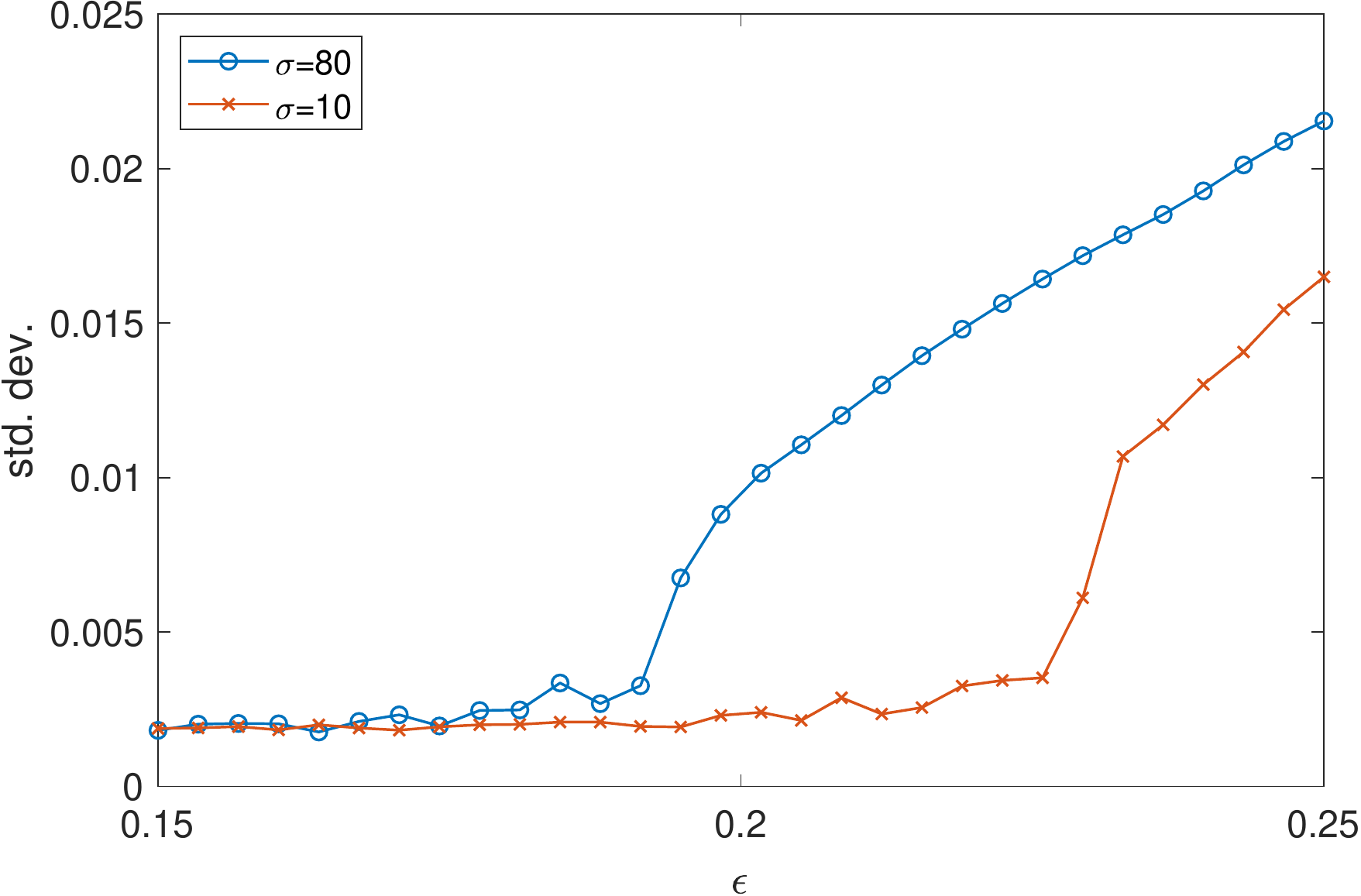}
\caption{Evidence of a supercritical Hopf bifurcation in a network of gap-junction coupled
Morris-Lecar neurons. On the vertical axis we plot the standard deviation of $\hat{s}$ 
(the mean of the $s_i$) over a period
of 10 seconds, having discarded the first 10 seconds as transients.}
\label{fig:hopfml}
\end{center}
\end{figure}

However, we can reproduce the movement of the Hopf bifurcation as $\sigma$ is varied 
in a network of Morris-Lecar neurons; see Fig.~\ref{fig:gapsmallg}.
We consider $N=2500$ and $\mk=100$ with a uniform degree
distribution on $[100-\sigma,100+\sigma]$. We obtain evidence of a supercritical
Hopf bifurcation as $\epsilon$ is increased as shown in Fig.~\ref{fig:hopfml}. 
On the vertical axis we plot the standard deviation of $\hat{s}$ over a period
of 10 seconds, having discarded the first 10 seconds as transients. As expected,
increasing $\sigma$ decreases the value of $g$ at which the bifurcation occurs.



\section{Gaussian or uniform distribution of $\eta_i$}
\label{sec:gauss}
All of the results so far have involved a Lorentzian distribution of a heterogeneous
parameter, either the input currents to theta neurons or to Morris-Lecar neurons.
In this section we investigate whether we obtain qualitatively similar results for
other distributions.

\begin{figure}
\begin{center}
\includegraphics[width=3.0in]{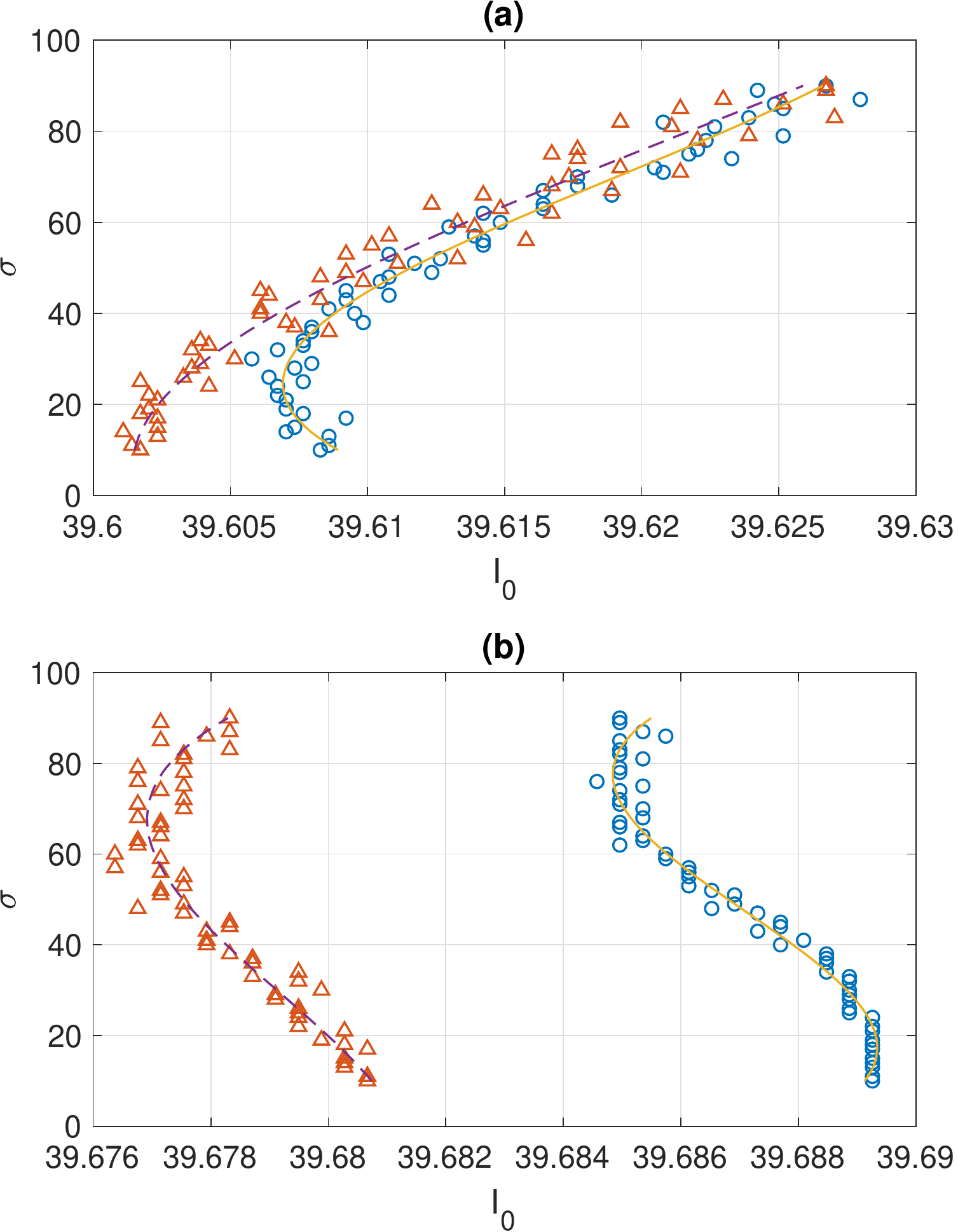}
\caption{SNIC bifurcations for gap-junction coupled Morris-Lecar neurons.
(a): $I_i$ chosen from a uniform distribution on $[-1/2,1/2]$ (blue circles, solid curve)
or a normal distribution with centre zero and standard deviation $1/3$ (red triangle, dashed
curve). (b): $I_i$ chosen from a uniform distribution on $[-1/8,1/8]$ (blue circles, solid curve)
or a normal distribution with centre zero and standard deviation $1/10$ (red triangle, dashed
curve). The curves result from fitting the values of $I_0$ as a cubic function of $\sigma$
and are to guide the eye. Parameters:
$N=2500,\epsilon=0.3$.}
\label{fig:gaunorm}
\end{center}
\end{figure}

Using either a Gaussian or uniform distribution of $I_i$ in a Morris-Lecar network
we obtained qualitatively the same results as in Figs.~\ref{fig:ML} and~\ref{fig:MLex}
for synaptic coupling (results not shown).
We investigated the effects shown in Fig.~\ref{fig:gaplargeg} in a network of
gap-junction coupled Morris-Lecar neurons for both uniform and normally
distributed $I_i$. The results are shown in Fig.~\ref{fig:gaunorm}. For each value of
$\sigma$ we created a network and a realisation of the $I_i$, and then used bisection in $I_0$
to approximately determine the transition from quiescence to periodic firing (with large period).
For broader distributions
(panel (a)) we obtained the same trend as in Fig.~\ref{fig:gaplargeg} while for narrower
distributions we seem to obtain the opposite trend (panel (b)). Note that
the theshold for firing for a single neuron is $I_0\approx 39.693455$, so all bifurcations occur
for $I_0$ less than this, in contrast with the results in Fig.~\ref{fig:gaplargeg}.
Such an effect has been seen before in excitable systems~\cite{lafcol10} indicating that the
Lorentzian distribution of heterogeneity, while providing analytical insight, may not
give generic results.

\begin{figure}
\begin{center}
\includegraphics[width=3.0in]{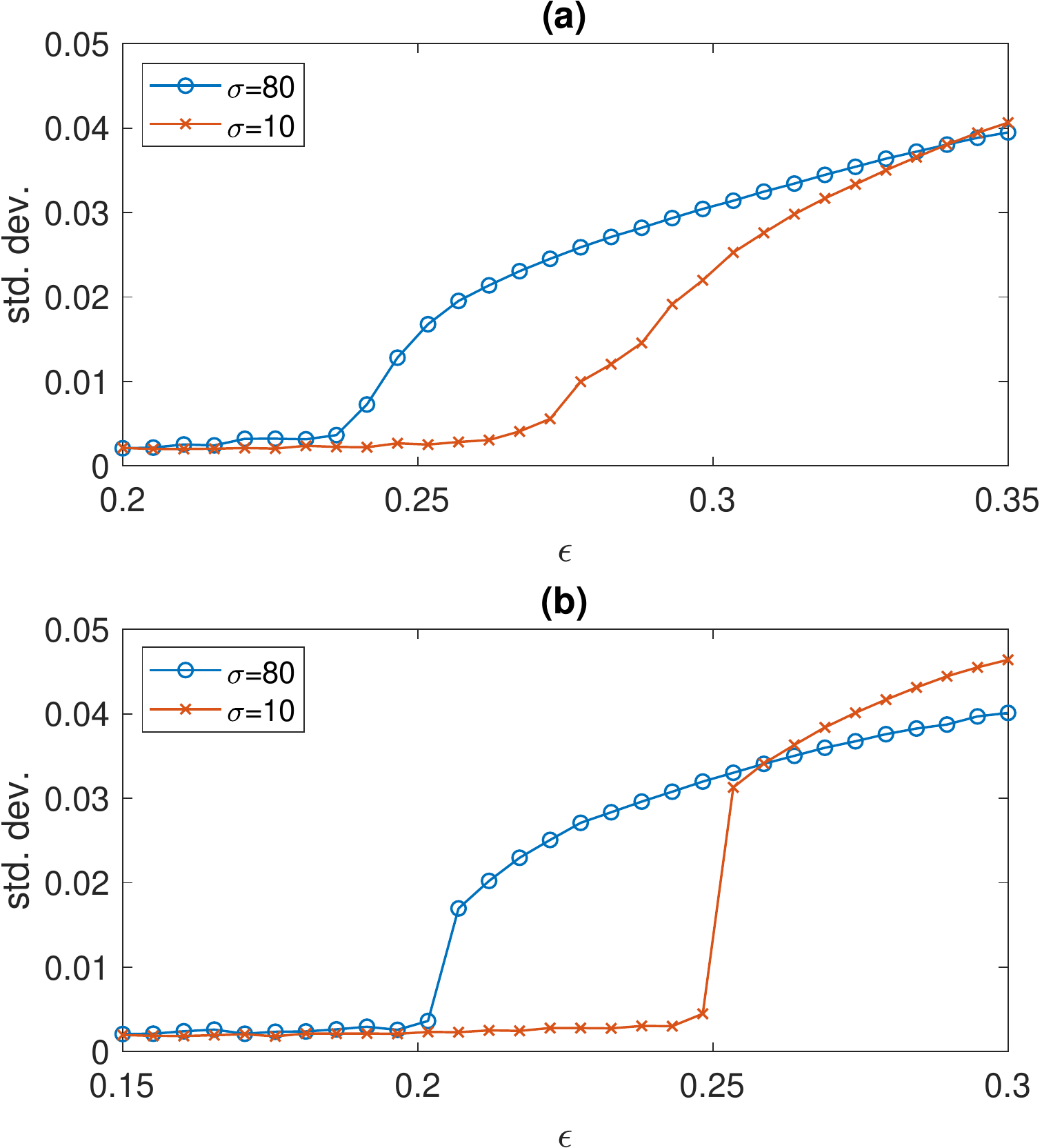}
\caption{Evidence of Hopf bifurcations in networks of gap-junction coupled
Morris-Lecar neurons. (a): $I_i$ taken from a unit Gaussian distribution. The Hopf bifurcations
appear to be supercritical.
(b): $I_i$ taken from a uniform distribution on $[-1,1]$. The Hopf bifurcation
for $\sigma=10$ appears to be subcritical, as explained in the text.
For each panel, on the vertical axis we plot the standard deviation of $\hat{s}$ 
(the mean of the $s_i$) over a period
of 10 seconds, having discarded the first 10 seconds as transients.
$N=2500$.}
\label{fig:gaphopf}
\end{center}
\end{figure}

We reproduced the results in Fig.~\ref{fig:gapsmallg} 
with the $I_i$ taken from a unit Gaussian (normal distribution), 
as shown in Fig.~\ref{fig:gaphopf}(a).
Choosing the $I_i$ from a uniform distribution 
on $[-1,1]$ we obtain the results in Fig.~\ref{fig:gaphopf}(b). Quasistatically sweeping
$\epsilon$ up and down for $\sigma=10$
we found a region of bistability between an approximate steady state and a 
macroscopic oscillation, suggesting that the Hopf bifurcation seen is subcritical.
For clarity, we only show the results of increasing $\epsilon$ for both networks. 
The effect of varying
the width of the in-degree distribution is the same: broadening the distribution moves the
Hopf bifurcation to a lower value of $\epsilon$.

\section{Summary}
\label{sec:summ}

We have derived approximate equations describing the expected dynamics of large networks
of theta neurons, under the assumption that the heterogeneous parameter has a Lorentzian
distribution. We have chosen the case of neutral assorativity within the networks
and independent in- and out-degrees, concentrating on the effects of varying the widths
of the in- and out-degree distributions. We have investigated both synaptic and gap
junctional coupling. Numerical bifurcation analysis has enabled us to determine the effects
of varying the degree distributions on the networks' dynamics.

For synaptically coupled inhibitory neurons, broadening the in-degree distribution
destroys macroscopic oscillations while for excitatory networks, it narrows
the range of values of mean input for which the network is bistable. The dynamics
are independent of the out-degree distribution. For gap junctional coupling,
broadening the degree distribution causes SNIC or Hopf bifurcations associated with
the onset of collective firing to move in parameter space. Most of the results are confirmed to
occur in networks of more realistic neurons with different distributions of heterogeneity,
but the use of Lorentzian distributions may not always give generic results.

This work could be generalised in several ways. One is to extend it to cover coupled
populations of excitatory and inhibitory neurons~\cite{rox11,coobyr19,lai17}. This 
would naturally result in more parameters to explore, as distributions of degrees relating
to up to four types of connection would need to be specified. Another area of interest
is the inclusion of noise in the neurons' dynamics. Noise is known to play a significant
role in neural dynamics~\cite{lailor09,ermter10} yet its presence would invalidate the use of the
Ott/Antonsen ansatz which lies behind the derivations in this paper. Goldobin et al.
have made progress in this area, developing perturbation theory away from noise-free
case~\cite{tyugol18,goldol20,goldiv21,goltyu18,goldol19}, and Ratas and 
Pyragas have applied these
ideas to networks of theta neurons~\cite{ratpyr19}.



{\bf Acknowledgement:} This work was partially supported by the Marsden Fund Council from Government funding, managed by Royal Society Te Aparangi, grant number 17-MAU-054.

\appendix

\section{Morris-Lecar equations}
\label{sec:app}
\subsection{Synaptic coupling}
For synaptic coupling the equations for $N$ neurons are
\begin{align}
   C\frac{dV_i}{dt} & = g_L(V_L-V_i)+g_{Ca}m_\infty(V_i)(V_{Ca}-V_i) \label{eq:dVdt} \\
 & +g_Kn_i(V_K-V_i)+I_0+I_i+\frac{\epsilon}{\mk}\sum_{j=1}^N A_{ij}s_j \nonumber \\
   \frac{dn_i}{dt} & = \frac{\lambda_0(w_\infty(V_i)-n_i)}{\tau_n(V_i)} \\
   \tau\frac{ds_i}{dt} & = s_\infty(V_i)-s_i 
\end{align}
where
\begin{align}
   m_\infty(V) & =0.5(1+\tanh{[(V-V_1)/V_2]}) \\
   w_\infty(V) & =0.5(1+\tanh{[(V-V_3)/V_4]}) \\
   \tau_n(V) & = \frac{1}{\cosh{[(V-V_3)/(2V_4)]}} \\
   s_\infty(V) & = 1+\tanh{(V/10)}.
\end{align}
Parameters are $V_1=-1.2,V_2=18,V_3=12,V_4=17.4,\lambda_0=1/15 msec^{-1},g_L=2,
g_K=8,
g_{Ca}=4,
V_L=-60,
V_{Ca}=120,
V_K=-80,
C=20\mu F/cm^2$ and these are unchanged throughout the paper. Voltages are in mV and conductances
are in $mS/cm^2$. These are taken from~\cite{tsukit06} but we have added synaptic dynamics.
The theshold for firing for single neuron is $I_0\approx 39.693455$.

For Fig.~\ref{fig:ML} we set $\epsilon=-1 mS/cm^2,I_0=41$. The $I_i$ were taking from a Lorentzian with mean
zero and HWHM $0.01$.
For Fig.~\ref{fig:MLex} we set $\epsilon=25 mS/cm^2,\tau=20,N=500$. 
The $I_i$ were taking from a Lorentzian with mean
zero and HWHM $0.01$.

\subsection{Gap junction coupling}

For gap-junctional coupling we use
\begin{align}
   C\frac{dV_i}{dt} & = g_L(V_L-V_i)+g_{Ca}m_\infty(V_i)(V_{Ca}-V_i)  \nonumber \\
 & +g_Kn_i(V_K-V_i)+I_0+I_i \nonumber \\
   & +\frac{\epsilon}{\mk}\sum_{j=1}^N A_{ij}(V_j-V_i) \\
   \frac{dn_i}{dt} & = \frac{\lambda_0(w_\infty(V_i)-n_i)}{\tau_n(V_i)} 
\end{align}
where $m_\infty(V),w_\infty(V)$ and $\tau_n(V)$ and all other parameters are as above.
For Fig.~\ref{fig:hopfml} we set $I_0=40$ and choose the $I_i$ from a Lorentzian with
mean zero and HWHM 0.5.

%


\end{document}